\documentclass[11pt,notitlepage,twoside]{article}
\usepackage{latexsym}
\usepackage{amsmath}
\usepackage{amsfonts}
\usepackage{amssymb}
\renewcommand{\a }{\alpha }
\renewcommand{\b }{\beta }
\renewcommand{\d}{\delta }

\newcommand{\D }{\Delta }

\newcommand{\e }{\varepsilon }

\renewcommand{\l }{\lambda }

\newcommand{\n }{\nabla }
\newcommand{\var }{\varphi }

\renewcommand{\th }{\theta }
\renewcommand{\o }{\omega }
\renewcommand{\O }{\Omega }

\newcommand{\ov}{\overline}

\newcommand{\be}{\begin{equation}}
\newcommand{\ee}{\end{equation}}
\newenvironment{pf}{\noindent{\bf Proof.}\enspace}{%\rule{2mm}{2mm}
\hfill$\Box$\medskip}
\newenvironment{pfn}[1]{\noindent{\bf Proof of {#1}\enspace}}{%\rule{2mm}{2mm}
\hfill$\Box$\medskip}
\newcommand{\R}{\mathbb{R}}
\newtheorem{thm}{Theorem}[section]
\newtheorem{pro}[thm]{Proposition}
\newtheorem{lem}[thm]{Lemma}
\newtheorem{rem}[thm]{Remark}

\numberwithin{equation}{section}

\author{Mohamed Ben Ayed$^a$\footnote{ Corresponding author. Fax : +216-74-274437, E-mail : \texttt{Mohamed.Benayed@fss.rnu.tn}.}\,\, \& Khalil El Mehdi$^{b,c}$\,\footnote{E-mail : \texttt{khalil@univ-nkc.mr}.} \\
{\footnotesize
a : D{\'e}partement de Math{\'e}matiques, Facult{\'e} des Sciences de Sfax, Route
Soukra, Sfax, Tunisia.}\\
{\footnotesize
b :  Facult\'e des Sciences et Techniques, Universit\'e de Nouakchott, BP 5026, Nouakchott, Mauritania.}\\
{\footnotesize
 c : The Abdus Salam ICTP, Mathematics Section,
 Strada Costiera 11, 34014 Trieste, Italy.}
}

\title { {\Large \textbf{On a Biharmonic Equation Involving Nearly Critical Exponent  }}\footnote{Work finished when the authors were visiting Mathematics Department of Roma University ``La Sapienza''. They would like to thank the Mathematics Department for its warm hospitality. The authors also thank Professors  M. Grossi and F. Pacella for their constant support.} }

\begin{document}

\date{ }

\maketitle
{\footnotesize

\noindent
{\bf Abstract.}
This paper is concerned with a biharmonic equation under the Navier boundary condition  $(P_{\mp\e})$ : $\D^2u=u^{\frac{n+4}{n-4}\mp\e}$, $u>0$ in $\O$ and $u=\D u=0$ on $\partial\O$, where $\O$ is a smooth bounded domain in $\R^n$, $n\geq 5$, and $\e >0$. We study the asymptotic behavior of solutions of $(P_{-\e})$ which are minimizing for the Sobolev quotient as $\e$ goes to zero. We show that such solutions concentrate around a point $x_0 \in \O$ as $\e\to 0$, moreover $x_0$ is a critical point of the Robin's function. Conversely, we show that for any nondegenerate critical point $x_0$ of the Robin's function, there exist solutions of $(P_{-\e})$ concentrating around $x_0$ as $\e\to 0$. Finally we prove that, in contrast with what happened in the subcritical equation $(P_{-\e})$, the supercritical problem  $(P_{+\e} )$ has no solutions which concentrate around a point of $\O$ as $\e\to 0$.

}

\section{Introduction and   Results }
\mbox{}
In this paper, we are concerned with the following semilinear biharmonic equation under the Navier boundary condition
$$
(P_{\mp\e} ) \quad \left\{
\begin{array}{cc}
 \D ^2 u = u^{p\mp\e},\,\, u>0 &\mbox{ in }\, \O \\
    \D u= u =0   \quad\quad     & \mbox{ on }\,  \partial  \Omega ,
\end{array}
\right.
$$
where $\O$ is a smooth bounded domain in $\R^n$, $n\geq 5$, $\e$ is a small positive parameter, and $p+1= 2n/(n-4)$ is the critical Sobolev exponent of the embedding  $H^2(\O )\cap H^1_0(\O ) \hookrightarrow L^{p+1}(\O )$.\\
When the biharmonic operator in $(P_{\mp\e})$ is replaced by the Laplacian operator, there are many works devoted to the study of the contrepart of $(P_{\mp\e})$, see for example \cite{AA}, \cite{AP}, \cite{BLR}, \cite{BEGR}, \cite{BP}, \cite{DFM1}, \cite{DFM2}, \cite{H}, \cite{L}, \cite{P1}, \cite{P2}, \cite{R2} and the references therein.\\
 When $\e\in (0,p)$, the mountain pass lemma proves the existence of solution to $(P_{-\e})$ for any domain $\O$. When $\e =0$, the situation is more complex, Van Der Vorst showed in \cite{V1} that if $\O$ is starshaped $(P_{0})$ has no solution whereas Ebobisse and Ould Ahmedou proved in \cite{EO} that  $(P_{0})$ has a solution provided that some homology group of $\O$ is nontrivial. This topological condition is sufficient, but not necessary, as examples of contractible domains $\O$ on which a solution exists show \cite{GGS}.\\
In view of this qualitative change in the situation when $\e=0$, it is interesting to study the asymptotic behavior of the subcritical solution $u_\e$ of  $(P_{-\e})$ as $\e\to 0$. Chou and Geng \cite{CG} made the first study, when $\O$ is a convex domain. The aim of the first result of this paper is to remove the convexity assumption on $\O$. To state this result, we need to introduce some notation.\\
Let us define on $\O$ the following Robin's function
$$
\var (x)=H(x,x),\quad \mbox{with}\quad H(x,y)=\frac{1}{|x-y|^{n-4}} -G(x,y),\,\,\mbox{for }\,\, (x,y)\in\O\times\O,
$$
where $G$ is the Green's function of $\D^2$, that is,
$$
\forall x\in\O \quad \left\{
\begin{array}{cc}
 \D ^2 G(x,.) = c_n\d_x &\mbox{ in }\, \O \\
    \D G(x,.)= G(x,.) =0      & \mbox{ on }\,  \partial  \Omega ,
\end{array}
\right.
$$
where $\d_x$ denotes the Dirac mass at $x$ and $c_n=(n-4)(n-2)|S^{n-1}|$.\\
Let
\begin{eqnarray}\label{e:11}
 \d _ {a,\l }(x) =  \frac {c_0\lambda
 ^{\frac{n-4}{2}}}{(1+\lambda^2|x-a|^2)^{\frac{n-4}{2}}},\,\,
 c_0=[(n-4)(n-2)n(n+2)]^{(n-4)/8},\,\,
 \l >0, \,\, a \in \R^n
\end{eqnarray}
It is well known (see \cite{Lin}) that $\d_{a,\l}$ are the only
solutions of
\begin{eqnarray*}
 \D^2 u =  u^{\frac{n+4}{n-4}},\quad  u>0 \mbox{  in  } \R^n, \quad
  \mbox{with } u\in L^{p+1}(\R^n) \quad \mbox{and } \D u \in L^2(\R^n)
\end{eqnarray*}
\noindent
and are also the only minimizers of the Sobolev inequality on the
whole space, that is
\begin{eqnarray}\label{e:12}
 S =\inf\{|\D u|^{2}_{L^2(\R^n)}|u|^{-2}_{L^{\frac{2n}{n-4}}(\R^n)}
,\, s.t.\, \D u\in L^2 ,u\in L^{\frac{2n}{n-4}} ,u\neq 0 \}.
\end{eqnarray}
\noindent
 We denote by  $P\d _{a,\l}$ the projection of the $\d
_{a,\l}$'s on $H^2(\O )\cap H^1_0(\O)$, defined by
$$ \D^2 P\d_{a,\l}=\D^2\d_{a,\l} \mbox{ in } \O  \mbox{ and }
\D P\d_{a,\l}=P\d_{a,\l}=0\mbox{ on } \partial \O.
$$
Let
\begin{align}
||u||&=\left(\int_\O |\D u|^2\right)^{1/2},\qquad u\in H^2(\O )\cap H^1_0(\O) \label{e:13}\\
(u,v)&=\int_\O \D u\D v,\qquad u,v \in H^2(\O )\cap H^1_0(\O) \label{e:14}\\
|u|_q&=|u|_{L^q(\O )} \label{e:15}.
\end{align}
Now we state the first result of this paper.
\begin{thm}\label{t:11} Assume that $n\geq 6$.
Let $(u_\e)$ be a solution of $(P_{-\e})$, and assume that
$$
{|| u_\e||^2}{|u_\e|_{p+1-\e}^{-2}}\to S \mbox{ as } \e \to
0,\leqno{(H)}
$$
where $S$ is the best Sobolev constant in $\R^n$ defined by \eqref{e:12}. Then (up to a
subsequence) there exist $a_\e\in \O$, $\l_\e > 0$, $\a_\e > 0$
and $v_\e$ such that $u_\e$ can be written as
$$
u_\e = \a_\e P\d _{a_\e , \l _\e}+v_\e
$$
with $\a_\e\to 1$, $||v_\e|| \to 0$, $a_\e \in
\O$ and $\l _\e d(a_\e, \partial\O ) \to +\infty$ as $\e \to 0$.\\
In addition, $a_\e$ converges to a critical point $x_0\in \O$ of
$\var$ and we have
$$
\lim_{\e\to 0}\e||u_\e||_{L^\infty (\O )}^2 = (c_1c_0^2/c_2)\var
(x_0),
$$
 where  $c_1=c_0^{2n/(n-4)} \int_{\R^n} \frac{d x
}{(1+|x|^2)^{(n+4)/2}}$, $c_2=(n-4)c_0^{2n/(n-4)}\int_{\R^n}
\frac{log(1+|x|^2) (1-|x|^2)}{(1+|x|^2)^{n+1}} dx $ and $c_0$ is
defined in \eqref{e:11}.
\end{thm}
\begin{rem}
It is important to point out that in the Laplacian case (see \cite{H}), the method of moving planes has been used to show that blowup points are away from the boundary of domain. The process is standard if domains are convex. For nonconvex regions, the method of moving planes still works in the Laplacian case through the applications of Kelvin transformations \cite{H}. For $(P_{-\e})$, the method of moving planes also works for convex domains \cite{CG}. However, for nonconvex domains, a Kelvin transformation does not work for  $(P_{-\e})$ because the Navier boundary condition is not invariant under the Kelvin transformation of biharmonic operator. Our method here is essential in overcoming the difficulty arising from the nonhomogeneity of Navier boundary condition under the Kelvin transformation.
\end{rem}
Our next result provides a kind of converse to Theorem \ref{t:11}.
\begin{thm}\label{t:12} Assume that $n\geq 6$, and $x_0\in\O$ is
a nondegenerate critical point of $\var$. Then
there exists an $\e_0 >0$ such that for each $\e\in (0,\e_0]$,
$(P_{-\e})$ has a solution of the form
$$
u_\e = \a_\e P\d _{a_\e , \l _\e}+v_\e
$$
with $\a_\e\to 1$, $||v_\e|| \to 0$, $a_\e \to
x_0$ and $\l _\e d(a_\e, \partial\O ) \to +\infty$ as $\e \to 0$.
\end{thm}
In view of the above results, a natural question arises:
are equivalent results true for slightly supercritical exponent?\\
The aim of the next result is to answer this question.
\begin{thm}\label{t:13}
Let $\O$ be any smooth bounded domain in $\R^n$, $n\geq 5$. Then
$(P_{+\e} )$ has no solution $u_\e$ of the form
$$
u_\e = \a_\e P\d _{a_\e , \l _\e}+v_\e
$$
with $||v_\e|| \to 0$, $\a_\e \to 1$, $a_\e \in
\O$ and $\l _\e d(a_\e,
\partial\O ) \to +\infty$ as $\e \to 0$.
\end{thm}
The proofs of our results are based on the same framework and
methods of  \cite{R1}, \cite{R2} and \cite{BEGR}. The next section
will be devoted to prove Theorem \ref{t:11}, while Theorems
\ref{t:12} and \ref{t:13} are proved in sections 3 and 4 respectively.

\section{Proof of Theorem \ref{t:11}}

Before starting the proof of Theorem \ref{t:11}, we need 
some preliminary results
\begin{pro}\label{p:21}\cite{BH}
Let $a\in \O$ and $\l>0$ such that $\l d(a,\partial \O)$ is large
enough. For $\th_{(a,\l)}=\d_{(a,\l)}-P\d_{(a,\l)}$, we have the
following estimates
$$
(a)\quad 0\leq \th_{(a,\l)}\leq\d_{(a,\l)}, \qquad (b) \quad
\th_{(a,\l)}=c_0\frac{H(a,.)}{\l^{\frac{n-4}{2}}}+f_{(a,\l)},$$ where
  $f_{(a,\l)}$ satisfies
$$f_{(a,\l)}=O\left(\frac{1}{\l^{\frac{n}{2}}d^{n-2}}\right),\quad
\l\frac{\partial f_{(a,\l)}}{\partial\l}=O\biggl(\frac{1}{
\l^{\frac{n}{2}}d^{n-2}}\biggr),\quad \frac{1}{\l}\frac{\partial
f_{(a,\l)}}{\partial a}=O\biggl(\frac{1}{
\l^{\frac{n+2}{2}}d^{n-1}}\biggr),$$ where $d$ is the distance
$d(a,\partial \O)$,
$$\mid\th_{(a,\l)}\mid_{L^{\frac{2n}{n-4}}}=O\bigl(\frac{1}{(\l
d)^{\frac{n-4}{2}}}\bigr), \quad
\mid\l\frac{\partial\th_{(a,\l)}}{\partial\l}\mid_{L^{\frac{2n}{n-4}}}=O\bigl(
\frac{1}{(\l d)^{\frac{n-4}{2}}}\bigr), \leqno{(c)}$$
 $$\mid\mid\th_{(a,\l)}\mid\mid=O\bigl(\frac{1}{(\l
d)^{\frac{n-4}{2}}}\bigr), \quad
 \mid\frac{1}{\l}\frac{\partial\th_{(a,\l)}}{\partial a}
\mid_{L^{\frac{2n}{n-4}}}=O\bigl(\frac{1}{(\l d)^{\frac{n-2}{2}}}\bigr).$$
\end{pro}

\begin{pro}\label{p:22}
Let $u_\e$ be a solution of $(P_{-\e})$ which satisfies $(H)$. Then,
we have
$$(a) \quad ||u_\e||^2\to S^{n/4}, \quad (b) \quad \int u_\e
^{p+1-\e}\to S^{n/4}.
$$
\end{pro}
\begin{pf}
Since $u_\e$ is a solution of $(P_{-\e})$, then we have
$||u_\e||^2=\int u_\e ^{p+1-\e} $. Thus, using the assumption
$(H)$, we derive that
$${|| u_\e||^2}{|u_\e|_{p+1-\e}^{-2}}=
||u_\e||^{\frac{2(p-1-\e)}{p+1-\e}} = S+o(1).
$$
Therefore $||u_\e||^2=\int u_\e ^{p+1-\e}=S^{n/4} +o(1)$. The
result follows.
\end{pf}
\begin{pro}\label{p:23}
Let $u_\e$ be a solution of $(P_{-\e})$ which satisfies $(H)$, and let $x_\e \in\O$ such that  $u_\e(x_\e)=|u_\e|_{L^\infty}:=M_\e$. Then, for $\e$
small, we have\\
$(a)\quad M_\e ^\e =1+o(1)$.\\
$(b) \quad u_\e$ can be written as
$$
u_\e = P\d _{x_\e , \tilde \l _\e}+\tilde v_\e,
$$
with $||\tilde v_\e|| \to 0$, where $\tilde \l _\e
=c_0^{2/(4-n)}M_\e ^{(p-1-\e)/4}$.
\end{pro}
\begin{pf}
First of all, we prove that $M_\e\to +\infty$ as $\e\to 0$. To this end, arguing by contradiction, we suppose that $M_{\e_n}$ remains bounded for a sequence $\e_n\to 0$ as $n\to +\infty$. Then, in view of elliptic regularity theory, we can extract a subsequence, still denoted by $u_{\e_n}$, which converges uniformly to a limit $u_o$. By Proposition \ref{p:22}, $u_0\neq 0$, hence by taking limit in $(H)$ we find that $u_0$ achieves a best Sobolev constant $S$, a contradiction to the fact that $S$ is never achieved on a bounded domain \cite{V2}.\\  
Now we define the rescaled functions
\begin{eqnarray}
\o _\e (y)=M_\e ^{-1}u_\e \left(x_\e + M_\e ^{(1+\e-p )/4}
y\right), \quad y\in \O _\e = M_\e ^{(p-1-\e )/{4}}(\O - x_\e ),
\end{eqnarray}
$\o _\e$ satisfies
\begin{eqnarray}
\left\{
\begin{array}{cccccc}
\D^2 \o _\e &=& \o _\e ^{p-\e},& 0< \o _\e \leq 1& \mbox{in}&\O _\e\\
\o _\e (0) &=&1,              & \D \o_\e =  \o _\e = 0    &
\mbox{on}&\partial \O _\e .
\end{array}
\right.
\end{eqnarray}
Following the same argument as in Lemma 2.3 \cite{BEH}, we have
$$
M_\e ^{(p-1-\e )/4}d(x_\e , \partial \O ) \to +\infty \quad
\mbox{as
  } \e \to 0.
$$
Then it follows from  standard elliptic theory that there exists a
positive function $\o $ such that (after passing to a subsequence)
$\o _\e \to \o $ in $C^4_{loc}(\R^n)$, and $\o$ satisfies
\begin{eqnarray*}
\left\{
\begin{array}{ccccc}
\D^2 \o &=&\o ^p,& 0\leq \o \leq 1 & \mbox{in }\R^n\\
\o (0)&=&1,&\n\o (0) =0. &
\end{array}
\right.
\end{eqnarray*}
It follows from \cite{Lin} that $\o$ writes as
$$
\o (y)= \d _{0,\a _n }(y), \quad \mbox{with } \a _n
=c_0^{2/(4-n)}.
$$
Observe that, for $y=M_\e ^{(p-1-\e)/4} (x-x_\e)$, we have
\begin{eqnarray}\label{*1}
M_\e \d_{0,\a_n}(y) =\frac{M_\e c_0 \a_n ^{(n-4)/2}}{\left(1+\a_n
^2 M_\e ^{(p-1-\e)/2}|x-x_\e|^2\right)^{(n-4)/2}}
 =M_\e ^{\e(n-4)/8}\d_{x_\e,\tilde \l _\e}
 \end{eqnarray}
  with  $\tilde \l _\e =\a_n M_\e ^{(p-1-\e)/4}$. Then,
$$
w_\e(y)-\d_{0,\a_n}(y)=M_\e ^{-1}\left( u_\e(x)-M_\e
^{\e(n-4)/8}\d_{x_\e,\tilde \l _\e}\right).
$$
Let us define
$$
u_\e ^1(x) =u_\e (x) - M_\e ^{\e(n-4)/8}P \d_{x_\e,\tilde \l
_\e}(x),
$$
we need to compute
$$
||u_\e ^1||^2=||u_\e ||^2+ M_\e ^{\e(n-4)/4}||P\d_{x_\e,\tilde \l
_\e}||^2-2M_\e ^{\e(n-4)/8}(u_\e, P\d_{x_\e,\tilde \l _\e}).
$$
On one hand, we have
\begin{align*}
(u_\e, P\d_{x_\e,\tilde \l _\e})& = \int_\O u_\e
(x)\d_{x_\e,\tilde \l_\e} ^{(n+4)/(n-4)}(x)\\
 & = \int_{\O_\e}
u_\e(x_\e+ M_\e ^{(1+\e-p)/4}y)\d_{x_\e,\tilde \l_\e}
^{(n+4)/(n-4)}(x_\e+
M_\e ^{(1+\e-p)/4}y)M_\e ^{n(1+\e-p)/4}dy\\
 & =\int_{\O_\e} M_\e ^{\e(n-4)/8}
w_\e(y)\d_{0,\a_n}^{(n+4)/(n-4)}(y)dy\\
&=\int_{B(0,R)} M_\e ^{\e(n-4)/8}
w_\e(y)\d_{0,\a_n}^{(n+4)/(n-4)}(y)dy + \int_{\O_\e\diagdown B(0,R)} M_\e ^{\e(n-4)/8}
w_\e(y)\d_{0,\a_n}^{(n+4)/(n-4)}(y)dy,
 \end{align*}
where $R$ is a large positive constant such that $\int_{\R^n\diagdown B(0,R)}\d_{0,\a_n}^{2n/(n-4)}=o(1)$.\\
Since
$$
\int_{\O_\e}w_\e^{2n/(n-4)}=M_\e^{-\e n/4}\int_\O u_\e^{2n/(n-4)}\leq c,
$$
using Holder's inequality we derive that
$$
\int_{\O_\e \diagdown B(0,R)}w_\e \d_{0,\a_n}^{(n+4)/(n-4)}=o(1).
$$
Now, since $w_\e \to \d_{0,\a_n}$ in $C^4_{loc}(\R^n)$, we obtain
$$
(u_\e, P\d_{x_\e,\tilde \l _\e})  =  M_\e ^{\e(n-4)/8} (S^{n/4}+o(1)).
$$
On the other hand, one can easy verify that
$$
||P\d_{x_\e,\tilde \l _\e})||^2=S^{n/4}+o(1).
$$
Thus
\begin{eqnarray}\label{*2}
||u_\e ^1||^2 = ||u_\e||^2-M_\e ^{\e(n-4)/4} (S^{n/4}+o(1)),
\end{eqnarray}
and, using the fact that $||u_\e ^1||^2\geq 0$ and Proposition
\ref{p:22}, we derive that
$$
M_\e ^{\e (n-4)/4} \leq 1+o(1).
$$
But, since $M_\e \to \infty$, we have $M_\e ^\e \geq 1$ and therefore claim $(a)$ follows.\\
Now we are going to prove claim $(b)$. Observe that,
using Proposition \ref{p:22} and claim $(a)$, \eqref{*2} becomes
$$
||u_\e ^1||^2 =(S^{n/4}+o(1)) - (S^{n/4}+o(1))=o(1).
$$
Thus claim $(b)$ follows.
\end{pf}
\begin{pro}\label{p:24}
Let $u_\e$ be a solution of $(P_{-\e})$ which satisfies $(H)$.
Then, there exist $a_\e\in \O$, $\a_\e > 0 $, $\l_\e >0$ and
$v_\e$ such that
$$
u_\e = \a_\e P\d _{a_\e , \l _\e}+ v_\e
$$
with $\a_\e \to 1$, $|a_\e -x_\e|\to 0$, $\l_\e d(a_\e,\partial
\O) \to \infty$, $\tilde \l _\e/\l_\e \to 1$ and $||v_\e||\to 0$.
Furthermore, $v_\e$ satisfies
$$
 (v,P\d_{a_\e,\l_\e})=(v,\partial P\d_{a_\e, \l_\e}/\partial
\l_\e)=0,\,  (v,\partial P\d_{a_\e, \l_\e}/\partial a)=0 .
\leqno{(V_0)}
$$
\end{pro}
\begin{pf}
By Proposition \ref{p:23}, $u_\e$ can be written as $u_\e =
P\d_{x_\e,\tilde \l _\e} +\tilde v_\e $ with $||\tilde v_\e ||\to
0$, $\tilde \l _\e d(x_\e,\partial \O)\to \infty$ as $\e \to 0$.
Thus, the following minimization problem
$$
\min\{ ||u_\e - \a P\d_{a,\l}||, \a >0 , a\in \O, \l > 0\}
$$
has a unique solution $(\a_\e, a_\e,\l_\e)$. Then, for $v_\e =
u_\e -\a_\e P\d_{a_\e,\l_\e}$, we have $v_\e$ satisfies $(V_0)$.
From the two forms of $u_\e$, one can easy verify that
$$ ||P\d_{x_\e,\tilde \l _\e}-P\d_{a_\e,\l_\e}||=o(1).
$$
Therefore, we derive that $|a_\e -x_\e|=o(1)$ and $\tilde \l _\e
/\l_\e =1+o(1))$. The result follows.
\end{pf}

Next, we state a result which its proof is similar to the proof of
Lemma 2.3 of \cite{BEGR}, so we will omit it.
\begin{lem}\label{l:25}
$\l_\e^\e =1+o(1)$ as $\e$ goes to zero implies that
\begin{eqnarray*}
\d_\e ^{-\e}-\frac{1}{c_0^\e\l_\e ^{\e(n-4)/2}}=O\left( \e
Log(1+\l_\e ^2 |x-a_\e|^2)\right) \quad \mbox{in}\quad \O.
\end{eqnarray*}
\end{lem}
We are now able to study the $v_\e$-part of $u_\e$ solution of $(P_{-\e})$.
\begin{pro}\label{p:26}
Let $(u_\e)$ be  a solution of $(P_{-\e})$ which satisfies $(H)$.
Then $v_\e$ occuring in Proposition \ref{p:24} satisfies
$$
|| v_\e ||\leq C \e + C\left( \frac{1}{(\l _\e
    d_\e )^{n-4}} \, (\mbox{if }n < 12)+
 \frac{1}{(\l _\e d_\e
    )^{\frac{n+4}{2}-\e(n-4)}}\, (\mbox{if } n \geq 12)\right),
$$
where $C$ is a positive constant independent of $\e$.
\end{pro}
\begin{pf}
Multiplying $(P_{-\e})$ by $v_\e$ and integrating on $\O$, we obtain
$$
 \int _{\O} \D u_\e .\D v_\e - \int _{\O} u_\e ^{p-\e} v_\e = 0.
$$
Thus
$$
\int _\O |\D v_\e |^2 - \int _\O \left( (\a _\e P\d _\e
)^{p-\e}+(p-\e)(\a _\e P\d _\e )^{p-1-\e}v_\e +O\left(\d _\e
^{p-2-\e}v_\e ^2 \chi _{|v_\e |<\d _\e} +|v_\e
|^{p-\e}\right)\right)v_\e=0.
$$
Using Lemma \ref{l:25}, we find
\begin{eqnarray}\label{e:27}
Q_\e (v_\e , v_\e ) - f_\e (v_\e ) + o(|| v_\e ||^2)=0,
\end{eqnarray}
with
$$
Q_\e (v,v)= || v||^2 - (p-\e)\int (\a _\e P\d _\e )^{p-1-\e
}v^2
$$
and
$$
f_\e (v)= \int (\a _\e P\d _\e )^{p-\e}v.
$$
We observe that
\begin{align*}
Q_\e &(v,v)= || v||^2 - p\int _\O (\a _\e P\d _\e )^{p-1-\e
}v^2 + O\left(\e || v||^2\right)\\
&= || v||^2 - p\a _\e ^{p-1-\e}\int _\O \left(\d _\e ^{p-1-\e}
+ O\left(\d _\e ^{p-2-\e}\th _\e\right)\right)v^2
+o\left(|| v||^2\right)\\
 & = ||v||^2 - \frac{p\a _\e ^{p-1-\e}}{c_0^\e \l _\e ^{\e
(n-4)/2}}\int _\O \d _\e ^{p-1}v^2  +O\left(\int _\O \left|\d _\e
^{-\e}-\frac{1}{c_0^\e \l _\e ^{\e (n-4)/2}}\right|\d _\e
^{p-1}|v|^2\right) +o\left(|| v||^2\right)
\end{align*}
Using Lemma \ref{l:25} and the fact that $\a_\e\to 1$, we find
$$
Q_\e (v,v) = Q_0 (v,v) + o\left(|| v||^2\right),
$$
with
$$
Q_0 (v,v)= || v||^2 - p\int _{\O} \d _\e ^{p-1}v^2.
$$
According to \cite{BE}, $Q_0$ is coercive, that is, there exists
some constant $c>0$ independent of $\e$, for $\e$ small enough,
such that
\begin{eqnarray}\label{e:28}
Q_0 (v,v) \geq c ||v||^2 \quad \forall v \in E_{(a _\e,\l _\e
)},
\end{eqnarray}
where
\begin{eqnarray}\label{es2}
 E_{(a _\e,\l _\e )}=\{v\in E/ v\,\,\mbox{satisfies}\,\, (V_0)\},
\end{eqnarray}
$(V_0)$ is the condition defined in Proposition \ref{p:24}.\\
We also observe that
\begin{align*}
f_\e (v) &= \a _\e ^{p-\e }\int _\O \left(\d _\e ^{p-\e}
+O\left(\d _\e ^{p-1-\e}\th _\e\right)\right)v\\
 & =\a _\e ^{p-\e }\left[\frac{1}{c_0^\e \l _\e ^{\e (n-4)/2}}\int _\O \d _\e
^p v + O\left( \e \int _\O Log (1+\l _\e ^2 |x-a_\e |^2 )\d _\e ^p
|v| + \int _\O\d _\e ^{p-1-\e}\th _\e |v| \right)\right].
\end{align*}
The last equality follows from Lemma \ref{l:25}. Therefore we can
write, with  $B = B(a_\e , d_\e )$
\begin{align*}
f_\e (v) & \leq c\left(\e || v|| + \int _{B}\d _\e
^{p-1-\e}\th _\e |v| + \int _{\R^n\diagdown B}\d _\e ^p
|v|\right)\\
 & \leq c|| v||\left(\e  + |\th _\e
|_{L^\infty}\left(\int _{B}\d _\e
^{(p-1-\e)\frac{2n}{n+4}}\right)^{\frac{n+4}{2n}} + \left(\int
_{\R^n\diagdown B}\d _\e
^{\frac{2n}{n-4}}\right)^{\frac{n+4}{2n}}\right).
\end{align*}
We notice that
\begin{eqnarray}\label{e:29}
\int _{\R^n\diagdown B}\d  _\e^{2n/(n-4)} = O\left(\frac{1}{(\l
_\e d_\e )^n}\right)
\end{eqnarray}
and
\begin{eqnarray}\label{e:210}
|\th_\e|_{L^\infty}\left(\int _{B}\d _\e
^{\frac{2n(p-1-\e)}{n+4}}\right)^{\frac{n+4}{2n}}\leq \frac{c}{(\l
_\e d_\e) ^{\frac{n+4}{2}-\e(n-4)}}(\mbox{ if } n \geq 12)+
\frac{c}{(\l _\e d_\e)^{n-4}} (\mbox{ if } n<12).
\end{eqnarray}
Thus we obtain
\begin{eqnarray}\label{e:211}
|f_\e (v)|\leq C || v||\left(\e + \frac{1}{(\l _\e d_\e
)^{n-4}} (\mbox{ if } n<12)+\frac{1}{(\l _\e d_\e
)^{\frac{n+4}{2}-\e(n-4)}} (\mbox{ if }n\geq 12)\right)
\end{eqnarray}
Combining \eqref{e:27}, \eqref{e:28} and \eqref{e:211}, we
obtain the desired estimate.
\end{pf}

Next we prove the following crucial result :
\begin{pro}\label{p:31}
 For $u_\e=\a_\e P\d_{a_\e,\l_\e}+v_\e$ solution of $(P_{-\e})$ with $\l_\e^\e =1 + o(1)$ as $\e$ goes to zero, we have the following estimate
$$
c_2\e+O(\e^2)-c_1\frac{H(a_\e,a_\e)}{\l_\e ^{n-4}}+o\left(\frac{1
}{(\l_\e d_\e)^{n-4}}\right)=0. \leqno{(a)}
$$
and for $n\geq 6$, we also have
$$
\frac{c_3}{\l_\e ^{n-3}}\frac{\partial H}{\partial a_\e}(a_\e,a_\e) +
O(\e ^2) + o\left(\frac{1}{(\l_\e d_\e)^{n-3}}\right)=0,
\leqno{(b)}
$$
where $c_1$, $c_2$ are the constants defined in Theorem \ref{t:11},
and where $c_3$ is a positive constant.
\end{pro}
\begin{pf}
We start by giving the proof of Claim $(a)$.
Multiplying the equation $(P_{-\e})$ by $\l_\e (\partial
P\d_\e)/(\partial\l_\e)$ and integrating on $\O$, we obtain
\begin{align}\label{e:31}
0&=\int _{\O}\D^2 u_\e \l _\e \frac{\partial P\d _\e}{\partial \l}
- \int _{\O}
u_\e ^{p-\e}\l _\e \frac{\partial P\d _\e}{\partial \l}\nonumber\\
&= \a _\e \int _{\O} \d _\e ^p \l _\e \frac{\partial P\d
_\e}{\partial \l}- \int _{\O} \left[ (\a _\e P\d _\e )^{p-\e} +
(p-\e ) (\a _\e P\d _\e )^{p-1-\e}v_\e \right.\nonumber\\
&\ \ \left.+O\left(\d _\e ^{p-2-\e}|v_\e |^2 +|v_\e
|^{p-\e}\chi_{\d_\e\leq |v_\e|}\right)\right]\l _\e \frac{\partial
P\d _\e}{\partial \l}.
\end{align}
We  estimate each term of the right-hand side in \eqref{e:31}.
First, using Proposition \ref{p:21}, we have
\begin{align*}
\int_{B^c} \d_\e ^p \l_\e\frac{\partial P\d_\e}{\partial \l_\e} &
\leq c \int_{B^c} \d_\e ^{p+1} = O(\frac{1}{(\l_\e
d_\e)^n})\\
\int _{B} \d _\e ^p \l _\e \frac{\partial P\d _\e}{\partial \l} &
= \int _{B} \d _\e ^p \l _\e \frac{\partial \d _\e}{\partial \l}+
\frac{(n-4)c_0}{2\l_\e ^{(n-4)/2}}\int _{B} \d _\e ^p H - \int_B
\d_\e ^p \l _\e \frac{\partial f_\e }{\partial \l},
\end{align*}
with $B=B(a_\e , d_\e )$.
Expanding $H(a_\e,.)$ around $a_\e$ and using Proposition
\ref{p:21}, we obtain
\begin{eqnarray*}
\int _{B} \d _\e ^p \l _\e \frac{\partial P\d _\e}{\partial \l}
=O\left(\frac{1}{(\l _\e d_\e )^n}\right)+\frac{(n-4)c_0}{2\l_\e
^{(n-4)/2}}H(a_\e,a_\e)\int_B\d_\e ^p +O\left(\frac{1}{(\l _\e
d_\e )^{n-2}}\right).
\end{eqnarray*}
 Therefore, estimating the integral, we obtain
\begin{eqnarray}\label{e:32}
\int _{\O} \d _\e ^p \l _\e \frac{\partial P\d _\e}{\partial \l}=
\frac{n-4}{2} c_1 \frac{H(a_\e ,a_\e )}{\l _\e ^{n-4}} +
O\left(\frac{1}{(\l _\e d_\e )^{n-2}}\right)
\end{eqnarray}
with ${c_1} = c_0^{2n/(n-4)}\int
_{\R^n}\frac{dx}{(1+|x|^2)^{(n+4)/2}}$.\\
Secondly, we compute
\begin{align}\label{e:33}
\int _{\O}(P\d _\e )^{p-\e}\l _\e \frac{\partial P\d _\e}{\partial
\l} & = \int _{\O}\left[\d _\e ^{p-\e}-(p-\e )\d _\e ^{p-1-\e}\th
_\e + O\left(\th _\e ^2 \d _\e ^{p-2-\e}+\th _\e
^{p-\e}\right)\right]\l _\e \frac{\partial P\d _\e
  }{\partial \l}\nonumber \\
&= \int _{B}\d _\e ^{p-\e}\l _\e \frac{\partial \d _\e}{\partial
\l} - \int _{B}\d _\e ^{p-\e}\l _\e \frac{\partial \th
_\e}{\partial \l}-(p-\e)\int _{B}\d _\e ^{p-1-\e}\th _\e \l
\frac{\partial \d _\e}{\partial \l}\\
&+ O\left(\int _{\O} \d _\e ^{p-1-\e}\th _\e |\l _\e
\frac{\partial \th _\e}{\partial \l}|+\int _{\O}\th _\e ^2 \d _\e
^{p-1-\e}+ \int _{\O}\th _\e ^{p-\e}\d _\e  + \frac{1}{(\l _\e
d_\e )^{n-\e\frac{n-4}{2}}}\right)\nonumber
\end{align}
and we have to estimate each term of the right hand-side of
\eqref{e:33}.\\
 Using the fact that $\l _\e \frac{\partial \d
_\e}{\partial \l} = \frac{n-4}{2} \left( \frac{1-\l _\e ^2 |x-a_\e
|^2}{1+\l _\e ^2 |x-a_\e |^2}\right) \d _\e $, we derive that
\begin{align}\label{e:34}
\int _{B}\d _\e ^{p-\e}\l _\e \frac{\partial \d _\e}{\partial
\l}&= \frac{n-4}{2}\frac{c_0^{p+1-\e}}{\l _\e ^{\frac{\e
(n-4)}{2}}}\int _{\R^n}\frac{1}{(1+|x|^2)^{n-\frac{\e (n-4)}{2}}}
\frac{1-|x|^2}{1+|x|^2}dx +O\left(\frac{1}{(\l _\e
    d_\e )^{n-\e(n-4)}}\right)\nonumber\\
 & = \frac{n-4}{2\l _\e ^{\e (n-4)/2}}\left(-c_2 \e + O(\e
^2 )\right)+O\left(\frac{1}{(\l _\e d_\e )^{n-\e(n-4)}}\right)
\end{align}
with $c_2 =\frac{n-4}{2}c_0^{\frac{2n}{n-4}}\int
_{\R^n}\frac{Log(1+|x|^2)}{(1+|x|^2)^n} \frac{|x|^2-1}{|x|^2+1}dx
>0$.\\
For the other terms in \eqref{e:33}, using Proposition \ref{p:21},
we have
\begin{align}\label{e:35}
\int _{B}\d _\e ^{p-\e}\l _\e \frac{\partial \th _\e}{\partial \l}
&  = \int _{B}\d _\e ^{p-\e}\l _\e \frac{\partial}{\partial
\l}\left(\frac{c_0H}{\l_\e ^{(n-4)/2}}-f_\e\right)\nonumber\\
&= -\frac{n-4}{2}{c_0^{p+1-\e}} \frac{H(a_\e ,a_\e )}{\l _\e
^{(n-4)/2}}\int _{B}\left(\frac{\l _\e}{1+\l _\e ^2|x-a_\e
|^2}\right)^{(p-\e ) \frac{(n-4)}{2}} +O\left(\frac{1}{
(\l _\e d_\e )^{n-2}}\right)\nonumber\\
&= -\frac{n-4}{2} \frac{H(a_\e ,a_\e )}{\l _\e ^{n-4}}\frac{1}{\l
_\e ^{\e (n-4)/2}} \int_{B(0,\l d)} \frac{c_0^{p+1-\e}}{(1
+|x|^2)^{\frac{n+4}{2}-\e\frac{n-4}{2}}}+ O\left(\frac{1}{ (\l _\e
d_\e )^{n-2}}\right)\nonumber\\
 & =-\frac{n-4}{2} c_1 \frac{H(a_\e ,a_\e )}{\l _\e ^{n-4}}\frac{1}{\l
_\e ^{\e (n-4)/2}}+O\left(\frac{\e}{ (\l _\e d_\e
)^{n-4}}+\frac{1}{ (\l _\e d_\e )^{n-2}}\right)
\end{align}
and
\begin{align*}
(p-\e)\int _{B}\d _\e ^{p-1-\e}\th _\e & \l _\e \frac{\partial \d
_\e}{\partial \l} =(p-\e)\int _{B}\d _\e ^{p-1-\e}\frac{c_0}{\l_\e
^{(n-4)/2}}H \l _\e \frac{\partial \d _\e}{\partial
\l}+O\left(\int
_{B}\d _\e ^{p-\e}f_\e\right)\\
 & =(p-\e)\frac{c_0}{\l_\e ^{(n-4)/2}}H(a_\e,a_\e) \int _{B}\d _\e
^{p-1-\e} \l _\e \frac{\partial \d _\e}{\partial
\l}+O\left(\frac{1}{(\l_\e d_\e)^{n-2}}\right)\\
 & =\frac{c_0(p-\e)}{\l_\e ^{\e\frac{n-4}{2}}}\frac{H(a_\e,a_\e)}{\l_\e
^{n-4}}\int_{B(0,\l_\e d_\e)}
\frac{(n-4)(1-|x|^2)}{2(1+|x|^2)^{\frac{n+6}{2}-\e
\frac{n-4}{2}}}+O\left(\frac{1}{(\l_\e d_\e)^{n-2}}\right)\\
 & =-\frac{n-4}{2}\frac{c_1}{\l_\e ^{\e(n-4)/2}}\frac{H(a_\e,a_\e)}{\l_\e
 ^{n-4}}+O\left(\frac{\e}{(\l_\e d_\e)^{n-4}}+\frac{1}{(\l_\e
 d_\e)^{n-2}}\right).
\end{align*}
 \eqref{e:33}, \eqref{e:34}, \eqref{e:35} and additional integral
 estimates of the same type provide us with the expansion
\begin{align}\label{e:36}
\int _{\O}(P\d _\e )^{p-\e}\l _\e \frac{\partial P\d _\e}{\partial
\l}= & \frac{n-4}{2\l _\e ^{\e (n-4)/2}}\left[-c_2\e +O(\e ^2) +2
c_1 \frac{H(a_\e ,a_\e )}{\l _\e ^{n-4}}\right.\\
 & \left.+ O\left(\frac{1 }{(\l _\e d_\e)^{n-2}}+(if\, n=5)
\frac{1}{(\l_\e d_\e)^2}\right)\right].\notag
\end{align}
We note that
\begin{align}\label{e:*1}
\int _{\O}(P\d _\e )^{p-1-\e}v_\e & \l _\e \frac{\partial P\d
_\e}{\partial \l} = \int _{\O} \left(\d _\e ^{p-1-\e} +O\left( \d _\e ^{p-2-\e}\th _\e \right)\right)v_\e \l _\e
\frac{\partial P\d _\e}{\partial
  \l}\notag\\
 & = \int _{\O}(\d _\e )^{p-1-\e}v_\e \l _\e \frac{\partial \d
_\e}{\partial
  \l}- \int _{\O}(\d _\e )^{p-1-\e}v_\e \l _\e \frac{\partial \th _\e}{\partial
  \l}+O\left(\int _{\O} \d _\e ^{p-1-\e}|v_\e |\th _\e \right)\notag\\
 & = \int _{\O}(\d _\e )^{p-1-\e}v_\e \l _\e \frac{\partial \d
_\e}{\partial \l}+O\left(\frac{|| v_\e ||}{(\l _\e d_\e ^2
)^{(n-4)/2}}\left(\int _{B}\d _\e ^{(p-1-\e)\frac{2n}{n+4}}
\right)^{(n+4)/2n}\right)\notag\\
&\ \ + O\left(|| v_\e ||\left(\int _{B^c}\d _\e
^{(p-\e)\frac{2n}{n+4}}\right)^{(n+4)/(2n)}\right).
\end{align}
Using \eqref{e:29} and \eqref{e:210}, we derive that
\begin{align}\label{e:e2}
\int _{\O}(P\d _\e )^{p-1-\e}v_\e \l _\e \frac{\partial P\d
_\e}{\partial \l} = & \int _{\O}\d _\e ^{p-1-\e}v_\e \l _\e
\frac{\partial \d _\e}{\partial \l}+ O\left(\frac{|| v_\e ||}{(\l
_\e d_\e )^{\frac{n+4}{2}-\e(n-4)}}\right) \\
 & + (\mbox{if }n<12)O\left(
\frac{||v_\e ||}{ (\l _\e d_\e) ^{n-4}}\right).\notag
\end{align}
We also have, using Lemma \ref{l:25}
\begin{align*}
\int _{\O}(\d _\e )^{p-1-\e}v_\e \l _\e \frac{\partial \d
_\e}{\partial \l} & = \frac{1}{c_0^\e\l _\e ^{\frac{\e
(n-2)}{2}}}\int _{\O}\d _\e ^{p-1}\l _\e v_\e \frac{\partial \d
_\e}{\partial \l} + \int _{\O}\left(\d _\e ^{-\e}-\frac{1}{c_0^\e
\l _\e ^{\frac{\e(n-2)}{2}}}\right)\d
_\e ^{p-1}v_\e \frac{\partial \d _\e}{\partial \l}\\
&=O\left(\e \int _{\O} Log(1+\l _\e ^2 |x-a_\e |^2 )\d _\e ^p
|v_\e | \right)=O(\e || v_\e||).
\end{align*}
Noticing that, in addition, $\l _\e \frac{\partial
P\d_\e}{\partial\l}=O(\d _\e )$ and
\begin{eqnarray}\label{e:38}
\int _{\O}\d _\e ^{p-1-\e}|v_\e |^2 = O\left(|| v_\e
||^2\right),\quad
 \int _{\d < |v_\e |}|v_\e |^{p-\e}\d _\e = O\left(||
v_\e||^{p+1-\e}\right).
 \end{eqnarray}
\eqref{e:32}, \eqref{e:36},..., \eqref{e:38}, Proposition
\ref{p:26}, Lemma \ref{l:25} and the fact that $\l_\e ^\e =1+O(\e
Log\l_\e)$ prove
claim (a) of our  Proposition.\\
Now, since the proof of Claim (b) is similar to that of Claim (a),
we only point out some necessary changes in the proof. We multiply
the equation $(P_{-\e})$ by  $(1/\l_\e)(\partial P\d_\e/\partial
a)$ and we integrate on $\O$, thus we obtain a similar equation as
\eqref{e:31}. As \eqref{e:32}, we have
\begin{align}\label{e:*2}
\int_{\O}\d_\e ^p \frac{1}{\l}\frac{\partial P\d_\e}{\partial a}=
& \int_B \d_\e ^p \frac{1}{\l}\frac{\partial \d_\e}{\partial
a}-\frac{c_0}{\l_\e ^{(n-2)/2}}\int_B \d_\e ^p \frac{\partial
H}{\partial a}+O\left(\frac{1}{(\l_\e d_\e)^{n-1}}\right)\notag\\
 & = -\frac{c_1}{2\l_\e ^{n-3}}\frac{\partial H}{\partial a}(a_\e,
 a_\e)+O\left(\frac{1}{(\l_\e d_\e)^{n-1}}\right).
 \end{align}
 Now, to obtain the similar result as \eqref{e:36}, we need to
estimate the following quantities
\begin{align*}
\int_B \d_\e ^{p-\e} \frac{1}{\l}\frac{\partial \d_\e}{\partial
a} & = 0\\
\int_B \d_\e ^{p-\e} \frac{1}{\l}\frac{\partial \theta_\e}
{\partial a} & = \int_B \d_\e ^{p-\e}\left(\frac{c_0}{\l_\e
^{(n-2)/2}} \frac{\partial H}{\partial a}-\frac{1}{\l
}\frac{\partial
f}{\partial a}\right)\\
 & = \frac{1}{\l_\e ^{\e\frac{n-4}{2}}\l_\e
^{n-3}}\frac{\partial H}{\partial a}(a_\e, a_\e) \left( c_1
+O(\e)\right)+O\left(\frac{1}{(\l_\e
d_\e)^{n-1}}\right)\\
\int\theta \frac{1}{\l}\frac{\partial} {\partial
a}(\d_\e)^{(p-\e)} & =(p-\e)D\theta (a_\e)\int\d_\e ^{p-1-\e}
\frac{1}{\l_\e}\frac{\partial \d_\e}{\partial
a}(x-a_\e)+O\left(\frac{1}{(\l_\e d_\e)^{n-1}}\right)\\
 & =\frac{1}{\l_\e ^{\e\frac{n-4}{2}}\l_\e
^{n-3}}\frac{\partial H}{\partial y}(a_\e, a_\e) \left( c_1
+O(\e)\right)+O\left(\frac{1}{(\l_\e d_\e)^{n-1}}\right).
\end{align*}
These estimates imply that
\begin{align}\label{e:*3}
\int_\O P\d_\e ^{p-\e} \frac{1}{\l}\frac{\partial P\d_\e}{\partial
a}= & \frac{- c_1}{\l_\e ^{\e\frac{n-4}{2}}\l_\e
^{n-3}}\frac{\partial H}{\partial a} (a_\e,
a_\e)+O\left(\frac{\e}{(\l_\e d_\e)^{n-3}}+\frac{1}{(\l_\e
d_\e)^{n-1}}\right)\\
 & + (if \, n=6) O\left( \frac{1}{(\l_\e
d_\e)^4}\right).\notag
\end{align}
Now, it remains to prove the similar result as \eqref{e:*1}.
Using the same arguments, we obtain
\begin{eqnarray}\label{e:e3}
 \int_{\O} P\d_\e ^{p-1-\e}v_\e \frac{1}{\l_\e}\frac{\partial P\d_\e}{\partial
a}=O\left(\e ^2+ \frac{1}{(\l_\e d_\e)^{n+4-\e(n-4)}}+(if\,  n<
12) \frac{1}{(\l_\e d_\e)^{2(n-4)}}\right).
\end{eqnarray}
The proof of claim (b) follows.
\end{pf}\\
We are now able to prove Theorem \ref{t:11}.\\
\begin{pfn}{Theorem \ref{t:11}}
Let $(u_\e)$ be a solution of $(P_{-\e})$ which satisfies $(H)$.
Then, using Proposition \ref{p:24}, $u_\e=\a_\e
P\d_{a_\e,\l_\e}+v_\e$ with $\a_\e\to 1$, $\l_\e d(a_\e, \partial
\O) \to \infty$, $v_\e$ satisfies $(V_0)$ and $|| v_\e||\to 0$. By
Propositions \ref{p:23} and \ref{p:24}, we have $\l_\e ^\e\to 1$.
Now, using claim $(a)$ of Proposition \ref{p:31}, we derive that
\begin{eqnarray}\label{e:39}
\e=\frac{c_1}{c_2}\frac{H(a_\e,a_\e)}{\l_\e
^{n-4}}+o(\frac{1}{(\l_\e d_\e)^{n-4}})=O(\frac{1}{(\l_\e
d_\e)^{n-4}}).
\end{eqnarray}
 Now claim (b) implies that
\begin{eqnarray}\label{e:310}
 \frac{\partial H}{\partial a_\e}(a_\e,a_\e)=o(\frac{1}{d_\e
 ^{n-3}}).
 \end{eqnarray}
Using \eqref{e:310} and the fact that for $a$ near the boundary $\partial H/\partial a (a_\e,a_\e)\sim cd^{3-n}(a,\partial\O)$, we derive that $a_\e$ is away from the boundary and it converges to a
critical point $x_0$ of $\var$.\\
 Finally, using \eqref{e:39}, we obtain
 $$
\e \l_\e ^{n-4} \to \frac{c_1}{c_2}\var (x_0) \mbox{ as } \e \to 0.
$$
Thus in order to complete the proof of our theorem, it only
remains to show that
 \begin{eqnarray}\label{e:311}
M_\e :=||u_\e||_{L^\infty}\sim c_0\l_\e^{(n-4)/2} \quad \mbox{as}
\quad \e\to 0.
\end{eqnarray}
Using Propositions \ref{p:23} and \ref{p:24}, we derive that $c_0^{2}||u_\e||_{L^\infty}
^{-2}\l_\e ^{n-4} \to 1$. Hence \eqref{e:311} follows. This
concludes the proof of Theorem \ref{t:11}.
\end{pfn}

\section{Proof of Theorem \ref{t:12}}

Let $x_0$ be a nondegenerate critical point of $\var$. It is easy
to see that $d(a, \partial \O) > d_0 > 0$ for $a$ near $x_0$. We will take a
function $u=\a P\d_{(a,\l)} +v$ where $(\a-\a_0)$ is very small,
$\l$ is large enough, $||v||$ is very small, $a$ is close to $x_0$
and $\a_0=S^{-n/8}$ and we will prove that we can choose the
variables $(\a, \l, a, v)$ so that $u$ is a critical point of
$J_\e$, where
$$
J_\e (u)=\left(\int_\O|\D u|^2\right) \left(\int_\O |u|^{p+1-\e}\right)^{-2/(p+1-\e)}
$$
is the functional corresponding to problem $(P_{-\e})$.\\
Let
\begin{align*}
M_\e =\{ (\a, \l, a, v)\in \R^*_+\times \R^*_+\times \O\times E/ & \ |\a-\a_0|< \nu_0,\  d_a > d_0,\  \l >
\nu_0^{-1}, \\
 & \e Log\l < \nu_0,\  ||v||<\nu_0 \mbox{ and } v\in
E_{(a,\l)}\},
\end{align*}
where $\nu_0$ and $d_0$ are two suitable positive constants and where $d_a=d(a,\partial\O )$.\\
 Let us define
the functional
\begin{eqnarray*}
K_\e :  M_\e \to \R, \quad  K_\e (\a, a, \l, v) = J_\e(\a P\d_{(a,
\l)} +v).
\end{eqnarray*}
Notice that $(\a, \l, a, v)$ is a critical point of $K_\e$
if and only if $u=\a P\d_{(a,\l)} + v$ is a critical point
of $J_\e$ on $E$. So this fact allows us to look for
critical points of $J_\e$ by successive optimizations with
respect to the different parameters on $M_\e$.\\
First, arguing as in Proposition 4 of \cite{R1} and using
computations performed in the previous sections, we observe
that the following problem
$$\min \{ J_\e(\a P\d_{(a, \l)} +v),
\ v \in E_{(a,\l)} \mbox{ and } ||v||<\nu_0\}
$$
is achieved by a unique function $\ov{v}$ which satisfies the
estimate of Proposition \ref{p:26}. This implies that there exist
$A$, $B$ and $C_i$'s such that
\begin{eqnarray}\label{ou1}
\frac{\partial K_\e}{\partial v}(\a, \l, a, \ov{v})=\n J_\e (\a
P\d_{(a, \l)} +\ov{v})=A  P\d_{(a, \l)}+B \frac{\partial
}{\partial \l}P\d_{(a, \l)}+\sum_{i=1}^nC_i \frac{\partial
}{\partial a_i}P\d_{(a, \l)},
\end{eqnarray}
where $a_i$ is the $i^{th}$ component of $a$.\\
 Now, we need to estimate the constants $A$, $B$, and $C_i$'s. For
this purpose, we take the scalar product of $\n J_\e (\a P\d_{(a,
\l)} +\ov{v})$ with $ P\d_{(a, \l)}$,  $ (\partial P\d_{(a,
\l)})/(\partial \l)$ and $ (\partial P\d_{(a,
\l)})/(\partial a_i)$ with $i=1,...,n$. Thus, we get a
quasi-diagonal system whose coefficients are given by
\begin{eqnarray*}
||P\d_{(a,\l)}||^2 =\ov{c}_1+O(\frac{1}{\l ^{n-4}}), \quad
\left(P\d_{(a, \l)},  \frac{\partial}{\partial \l}P\d_{(a,
\l)}\right)= O(\frac{1}{\l ^{n-3}}),
\end{eqnarray*}
$$ \left(P\d_{(a, \l)},  \frac{\partial}{\partial
a_i}P\d_{(a, \l)}\right) = O(\frac{1}{\l ^{n-4}}),\quad
\bigg|\bigg| \frac{\partial}{\partial \l}P\d_{(a,
\l)}\bigg|\bigg|^2 = \frac{ \ov{c}_2}{\l^2}+O(\frac{1}{\l ^{n-2}}),$$
$$ \left(\frac{
\partial}{\partial \l} P\d_{(a, \l)},
\frac{\partial}{\partial a_i}P\d_{(a, \l)}\right)=O(\frac{1}{\l
^{n-3}}), \quad \left( \frac{\partial}{\partial
a_j}P\d_{(a, \l)}, \frac{\partial }{\partial
a_i}P\d_{(a, \l)}\right)=\ov{c}_3\l^2\d_{ij} + O(\frac{1}{\l ^{n-5}})
$$
where $\ov{c}_i$'s are positive constants and $\d_{ij}$ is the
Kronecker symbol.\\
The other hand-side is given by $(\n J_\e(\a P\d_{(a,\l)}
+\ov{v}), \psi) $ where $\psi=P\d_{(a, \l)}$,  $(\partial
P\d_{(a, \l)})/(\partial \l)$, $ (\partial P\d_{(a,
\l)})/(\partial a_i)$ with $i=1,...,n$. \\
Observe that, for $u=\a P\d_{(a, \l)}+\ov{v}$,
\begin{eqnarray*}
(\n J_\e(u), \psi) = 2J_\e (u) \left( \a (P\d_{(a, \l)},
\psi)-J_\e (u)^{(p+1-\e)/2}\int _\O u^{p-\e}\psi\right).
\end{eqnarray*}
Expanding $J_\e$, we obtain
\begin{eqnarray}\label{e:*..2}
J_\e(\a P\d_{(a,\l)}+\ov{v})=S+O\left( \e+\e Log\l +\frac{1}{\l
^{n-4}}\right).
\end{eqnarray}
Now, using \eqref{e:211}, \eqref{e:32}, \eqref{e:36},
\eqref{e:e2}, \eqref{e:*2}, \eqref{e:*3}, \eqref{e:e3} and
Proposition \ref{p:26}, we derive that, after taking the following
change of variable: $\a =\a_0 +\b$,
\begin{align*}
(\n J_\e(u), P\d) & =O\left(\e Log \l +|\b|+\frac{1}{\l ^{n-4}}\right)\\
(\n J_\e(u), \partial P\d/\partial \l) & =O\left(\frac{\e}{\l}+\frac{1}{\l ^{n-3}}\right)\\
(\n J_\e(u),  \partial P\d/\partial a_j) & =O\left(\l\e
^2+\frac{1}{\l ^{n-4}}\right), \mbox{ for each } j=1,...,n.
\end{align*}
The solution of the system in $A$, $B$, and $C_i$'s shows that
\begin{eqnarray*}
A =O\left(\e Log \l +|\b|+\frac{1}{\l ^{n-4}}\right), \quad B
=O\left(\l \e +\frac{1}{\l ^{n-5}}\right),\quad C_j
=O\left(\frac{\e ^2}{ \l} +\frac{1}{\l ^{n-2}}\right).
\end{eqnarray*}
Now, to find critical points of $K_\e$, we have to solve the following system
$$
(E_1)\qquad
\begin{cases}
\frac{\partial K_\e}{\partial \a} + \left(\frac{\partial K}{\partial v}, \frac{\partial \bar{v}}{\partial\a}\right)&  =0\\
\frac{\partial K_\e}{\partial \l}+ \left(\frac{\partial K}{\partial v}, \frac{\partial \bar{v}}{\partial\l}\right) & = 0\\
\frac{\partial K_\e}{\partial a_j}+ \left(\frac{\partial K}{\partial v}, \frac{\partial \bar{v}}{\partial a_j}\right) & = 0,
\mbox{ for  } j=1,...,n.
\end{cases}
$$
Taking the derivatives, with respect to the different parameters on $M_\e$, of the following equalities
\begin{eqnarray*}
(V_0)\quad (\bar{v},P\d_{a_\e,\l_\e})=(\bar{v},\partial P\d_{a_\e,
\l_\e}/\partial \l_i)=0,\,  (\bar{v},\partial P\d_{a_\e,
\l_\e}/\partial a_i)=0 \mbox{ for }i=1,...,n
\end{eqnarray*}
and using \eqref{ou1}, we see that system $(E_1)$ is equivalent to
$$
(E_2)\qquad
\begin{cases}
\frac{\partial K_\e}{\partial \a} & =0\\
\frac{\partial K_\e}{\partial \l} & = B\left(\frac{\partial ^2
P\d}{\partial \l ^2}, \bar{v}\right)+\sum_{i=1}^n C_i
\left(\frac{\partial ^2 P\d}{\partial \l \partial a_i},
\bar{v}\right)\\
\frac{\partial K_\e}{\partial a_j} & =B\left(\frac{\partial ^2
P\d}{\partial \l \partial a_j}, \bar{v}\right)+\sum_{i=1}^n C_i
\left(\frac{\partial ^2 P\d}{\partial a_i\partial a_j}, \bar{v}\right),
\mbox{ for each } j=1,...,n.
\end{cases}
$$
The same computation as in the proof of Proposition \ref{p:26}
shows that
\begin{align*}
\frac{\partial K_\e}{\partial \a} & =(\n J_\e(\a P\d +\ov{v}), P\d)\\
 & = 2J_\e(u) \left( \a S^{n/4} \left(1-\a ^{p-1} S^{ n/(n-4)}\right)
 + O\left(\e Log \l +\frac{1}{\l ^{n-4}}\right)\right).
 \end{align*}
Furthermore, using the estimates provided in the proof of claim
(a) of Proposition \ref{p:31}, we derive that
\begin{align*}
\l \frac{\partial K_\e}{\partial \l} & =\left(\n J_\e(\a P\d
+\ov{v}),
\l \frac{\partial P\d}{\partial \l}\right)\\
 & = (n-4)J_\e(u) \left( \a c_1\frac{H(a,a)}{\l ^{n-4}}\left(1-2
\a ^{p-1} S^{ n/(n-4)}\right) +c_2S^{\frac{n}{n-4}}\a ^p \e\right. \\
 & \left.+O\left(\e ^2 Log \l +\frac{\e Log \l}{\l ^{n-4}}+\frac{1}{\l
^{n-2}} \right)\right).
 \end{align*}
Following also the proof of claim (b) of Proposition \ref{p:31},
we obtain, for each $j=1,...,n$,
\begin{align*}
\frac{\partial K_\e}{\partial a_j} & =\left(\n J_\e(\a P\d +\ov{v}),
 \frac{\partial P\d}{\partial a_j}\right)\\
 & = - \frac{c\a }{\l ^{n-4}}\frac{\partial H}{\partial a}(a,a)\left(1-2
\a ^{p-1} S^{ n/(n-4)}\right) + O\left(\l \e ^2  +\frac{\e Log
\l}{\l ^{n-4}}+\frac{1}{\l ^{n-2}}+(\mbox{if }n=6)\frac{1}{\l ^3}
\right).
 \end{align*}
On the other hand, one can easy verify that
\begin{eqnarray}\label{e:R}
(i)\,\, ||\frac{\partial^2 P\d}{\partial \l^2}||=O\left(\frac{1}{\l^2}\right), \qquad (ii)\,\, ||\frac{\partial^2 P\d}{\partial \l \partial a_i}||=O(1), \qquad (iii)\,\, ||\frac{\partial^2 P\d}{\partial a_i \partial a_j}||=O(\l^2).
\end{eqnarray}
Now, we take the following change of variables: \\
\begin{eqnarray*}
\a= \a_0 +\b, \quad a=x_0+\xi, \quad \frac{1}{\l ^{\frac{n-4}{2}}}
= \sqrt{\frac{c_2}{c_1}}\left( \frac{1}{H(x_0,x_0)}+\rho\right)
\sqrt{\e}.
\end{eqnarray*}
Then, using estimates \eqref{e:R}, Proposition \ref{p:26} and the fact that $x_0$ is a nondegenerate critical point
of $\var$, the system $(E_2)$ becomes
$$
(E_3) \qquad
\begin{cases}
\b & = O\left( \e |Log \e|+|\b|^2\right)\\
\rho & =O\left( \e ^{2/(n-4)}+|\b|^2+|\xi|^2+\rho ^2\right)\\
\xi & =O\left( \e ^{2/(n-4)}+|\b|^2+|\xi|^2+ \rho
^2+(\mbox{if } n=6)\e ^{1/2}\right).
\end{cases}
$$
Thus Brower's fixed point theorem shows that the system $(E_3)$ has a solution $(\b^\e,
\rho^\e, \xi^\e)$ for $\e$ small enough such that
\begin{eqnarray*}
\b^\e = O(\e |Log \e|), \quad \rho^\e =O(\e^{2/(n-4)} +(\mbox{if
}n=6) \e ^{1/2}), \quad \xi^\e=O(\e^{2/(n-4)} +(\mbox{if }n=6) \e
^{1/2}).
\end{eqnarray*}
 By construction, the corresponding $u_\e$ is a critical point of
 $J_\e$ that is $w_\e = J_\e(u_\e)^{n/8}u_\e$ satisfies
 \begin{eqnarray}\label{e:ew}
 \D ^2 w_\e = |w_\e|^{8/(n-4) - \e}w_\e \mbox{ in }\O, \quad w_\e=\D
 w_\e=0 \mbox{ on } \partial \O.
 \end{eqnarray}
with $|w_\e^-|_{L^{2n/(n-4)}(\O)}$ very small, where $w_\e
^-=\max(0,-w_\e)$.\\
As in Propostion 4.1 of \cite{BEH2}, we prove that $w_\e^-=0$. Thus, since $w_\e$ is
a non-negative function which satisfies \eqref{e:ew}, the strong maximum
principle ensures that $w_\e > 0$ on $\O$ and then $u_\e$ is a
solution of $(P_{-\e})$, which blows-up at $x_0$ as $\e$ goes to
zero. This ends the proof of Theorem \ref{t:12}.

\section{Proof of Theorem \ref{t:13}}
First of all, we
can easily show that for
$u_\e$ satisfying the assumption of the theorem, there is a unique
way to choose $a_\e$, $\l _\e$ and $v_\e$ such that
\begin{eqnarray}\label{e:51}
u_\e = \a _\e P\d _{a_\e , \l _\e}+v_\e
\end{eqnarray}
with
\begin{eqnarray}\label{e:52}
\left\{
\begin{array}{lll}
\a _\e \in \R, \quad \a _\e \to 1\\
a_\e \in \O, \quad \l _\e \in \R^*_+,\quad \l _\e d(a_\e , \partial\O
)\to +\infty\\
v_\e \to 0 \quad\mbox{in } E:=H^2\cap H_0^1(\O), \quad v_\e \in
E_{a_\e ,\l _\e }
\end{array}
\right.
\end{eqnarray}
and
for any $(a,\l ) \in \O
\times \R^*_+$, $E_{(a,\l)}$ denotes the subspace of $E$
  defined by \eqref{es2}.\\
 In the following,
we always assume that $u_\e$, satisfying the assumption of
Theorem \ref{t:13}, is written as in \eqref{e:51}.
To simplify the notations, we set
 $\d _{a_\e , \l _\e }=\d _\e$, $P\d
_{a_\e , \l _\e }=P\d _\e$ and $\th _{a_\e ,\l _\e}=\th _\e$.\\
Now we are going to  estimate the $v_\e$ occurring in
\eqref{e:51}.
\begin{lem}\label{l:51}
Let $u_\e$ satisfying the assumption of Theorem \ref{t:13}. Then we have
$$
(i)\, \int _{\O}|\D u_\e |^2 \to S^{n/4}~; \quad (ii) \, \int _{\O}u_\e
^{p+1+\e}\to S^{n/4}
$$
as $\e \to 0$, $S$ denoting the Sobolev constant defined by \eqref{e:12}.
\end{lem}
\begin{pf}
We have
\begin{align*}
\int _{\O}|\D u_\e |^2 &= \int _{\O}|\D (\a _\e P\d _\e + v_\e )|^2\\
       &= \a _\e ^2 \int _{\O}|\D P\d _\e |^2 + \int _{\O}|\D v_\e |^2
       \quad \mbox{since } v_\e \in E_{a_\e ,\l _\e}.
\end{align*}
From the fact that $\d _\e$ satisfies $\D^2 \d _\e =\d ^p_\e$ in $\R^n$
and is a minimizer for $S$, we deduce that
$$
\int _{\R^n}|\D \d _\e |^2 = S^{n/4}.
$$
On the other hand, an explicit computation provides us with
$$
\int _\O |\D\d _{a,\l}|^2=\int _{\R^n}|\D\d _{a,\l}|^2 +
O\left(\frac{1}{(\l d(a,\partial\O ))^n}\right) \quad \mbox{as } \l
d(a,\partial\O )\to +\infty .
$$
Using Proposition \ref{p:21}, claim (i) is a consequence of
\eqref{e:52}. Claim (ii) follows from the fact that $u_\e$ solves
$(P_{+\e} )$.
\end{pf}\\
\begin{lem}\label{l:52}
Let $u_\e$ satisfying the assumption of Theorem \ref{t:13}. Then  $\l _\e$ occuring
in \eqref{e:51} satisfies
$$
\l _\e ^\e  \to 1  \quad \mbox{as } \e \to 0.
$$
\end{lem}
\begin{pf}
According to Lemma \ref{l:51}, we have
\begin{eqnarray}\label{e:53}
\int _{\O}u_\e ^{p+1+\e} = S^{n/4} + o(1) \quad \mbox{as } \e \to 0
\end{eqnarray}
and
\begin{align*}
\int _{\O}u_\e ^{p+1+\e}&= \int _{\O}(\a _\e P\d _\e + v_\e )^{p+\e}\a
_\e P\d _\e  + \int _{\O}u_\e ^{p+\e}v_\e \\
&=\a _\e ^{p+\e +1}\int _{\O}P\d _\e ^{p+\e +1}+
 \int _{\O} \D^2 u_\e v_\e + O\left(\int _{\O}P\d _\e
^{p+\e}|v_\e | + \int _{\O}\mid v_\e \mid ^{p+\e}P\d _\e \right) \\
&=\a _\e ^{p+\e +1}\int _{\O}P\d _\e ^{p+\e +1} + O\left(\l _\e
  ^{\frac{\e(n-4)}{2}}\int _{\O} P\d _\e ^p |v_\e |+ \l _\e
  ^{\frac{\e(n-4)}{2}}\int _{\O}\mid v_\e\mid ^{p+\e}P\d _\e ^{1-\e} +||v_\e
  ||\right)\\
&= \a _\e ^{p+\e +1}\int _{\O}P\d _\e ^{p+\e +1} + O\left(\l _\e
  ^{\e(n-4)/2}|v_\e |_{L^{p+1}} + \l _\e
  ^{\e(n-4)/2}|v_\e |^{p+\e}_{L^{p+1}} + ||v_\e ||\right).
\end{align*}
Thus
\begin{eqnarray}\label{e:54}
\int _{\O}u_\e ^{p+1+\e}= \a _\e ^{p+\e +1}\int _{\O}P\d _\e ^{p+\e +1} +
o\left(\l _\e ^{\e(n-4)/2} + 1\right).
\end{eqnarray}
We observe that
\begin{align}
\int_{\O}P\d _\e ^{p+1+\e} & =\int_{\O}(\d _\e - \th _\e )
^{p+1+\e} =\int _{\O}\d _\e ^{p+1+\e}+ O\left(\int _{\O}\d _\e
^{p+\e}\th _\e
 \right)\notag \\
  & =c_0^{p+\e+1}\int _{\R^n} \left(\frac{\l _\e}{1+\l _\e ^2|x-a_\e
|^2}\right)^{n+\e(n-4)/2}\notag\\
& +O\left(|\th _\e |_{L^\infty} \int _\O \left(\frac{\l _\e}{1+\l
_\e ^2|x-a_\e|^2} \right)^{\frac{(p+\e)(n-4)}{2}} +\frac{\l _\e
 ^{\frac{\e(n-4)}{2}}}{(\l _\e d_\e )^{n}}\right).\notag
\end{align}
 Using Proposition \ref{p:21}, we obtain
\begin{eqnarray*}
\int _{\O}P\d_\e^{p+1+\e}=  c_0^{p+1+\e} \l _\e
^{{\e(n-4)/2}}  \int _{\R^n} \frac{dx}{(1+|x|^2)^{n+\e(n-4)/2}}
+O\left(\frac{\l _\e ^{\frac{\e(n-4)}{2}}}{(\l _\e d_\e
)^{n-4}}\right).
\end{eqnarray*}
We note that
\begin{align*}
c_0^{p+1+\e}\int _{\R^n}\frac{dx}{(1+|x|^2)^{n+\e(n-4)/2}} &=
c_0^{p+1}\int _{\R^n}\frac{dx}{(1+|x|^2)^n} +O(\e)\\
&= S^{n/4} +O(\e).
\end{align*}
Therefore
\begin{eqnarray}\label{e:55}
\int _{\O}P\d _\e ^{p+1+\e}=\l _\e
^{\e(n-4)/2}\left(S^{n/4}+O(\e) + o(1)\right).
\end{eqnarray}
and \eqref{e:54} and \eqref{e:55} provide us with
\begin{eqnarray}\label{e:56}
\int _{\O}u_\e ^{p+1+\e}=\a _\e ^{p+1+\e} \l _\e ^{\e(n-4)/2}\left(S^{n/4}+o(1)\right) +o(1).
\end{eqnarray}
Combination of \eqref{e:53} and \eqref{e:56} proves the lemma.
\end{pf}\\

Next, as in Lemma 2.3 of \cite{BEGR}, we can easily prove the following estimate :
\begin{lem}\label{l:53} $\l_\e^\e=1+o(1)$ as $\e$ goes to zero implies that
$$
\d  _\e^\e (x) - c_0^\e \l _\e ^{\e (n-4)/2} = O\left(\e Log (1+\l _\e ^2
  |x-a_\e|^2)\right) \qquad \mbox{in } \O .
$$
\end{lem}
We are now able to study the $v_\e$-part of $u_\e$.
\begin{lem}\label{l:54}
Let $u_\e$ satisfying the assumption of Theorem \ref{t:13}. Then
$v_\e $ occuring in  \eqref{e:51} satisfies
$$
\int _{\O}|v_\e |^{p+1+\e} = o(1) \quad \mbox{as } \e \to
0.
$$
\end{lem}
\begin{pf}
We observe that
\begin{eqnarray*}
\int _{\O}u_\e ^{p+1+\e} = \int _{\O}(\a _\e P\d _\e
)^{p+1+\e} +  \int _{\O}|v_\e
|^{p+1+\e} + O\left( \int _{\O}(\a _\e P\d _\e
)^{p+\e}|v_\e | +   \int _{\O}|v_\e
|^{p+\e}\a _\e P\d _\e \right).
\end{eqnarray*}
We are going to estimate each term of the right hand-side in the
above equality.
\begin{align*}
\int _{\O}(\a _\e P\d _\e )^{p+1+\e}&= \a _\e ^{p+1+\e}\left[\int
  _{\O}\d _\e ^{p+1+\e}+ O\left( \int _{\O}\d _\e
^{p+\e}\th _\e \right)\right]\\
&=\a _\e ^{p+1+\e}\int _{\O}\d _\e ^{p+1+\e
      }+ o\left(\l _\e ^{\e (n-4)/2}\right)\\
 \int _{\O}(\a _\e P\d _\e )^{p+\e}|v_\e |&\leq \l _\e
 ^{\frac{\e (n-4)}{2}}||v_\e || = o(1)\\
\int _{\O}|v_\e |^{p+\e}\a _\e P\d _\e &\leq  \l _\e
 ^{\frac{\e (n-4)}{2}}||v_\e ||^{p+\e} = o(1)
\end{align*}
using Proposition \ref{p:21}, Holder inequality and Sobolev embedding theorem.
From Lemma \ref{l:51} we derive that
\begin{eqnarray}\label{e:57}
S^{n/4} +o(1) = (1+o(1))\int _{\O}\d  _\e^{p+1+\e}+  \int _{\O}|v_\e
|^{p+1+\e}.
\end{eqnarray}
As we have also
$$
\int _{\O}\d  _\e^{p+1+\e}= \l _\e ^{\frac{\e
    (n-4)}{2}}c_0^{\e}\int _{\O}\d _\e ^{p+1} + \int
_{\O}\left(\d _\e ^{p+1+\e}- c_0^{\e} \l _\e ^{\frac{\e
    (n-4)}{2}}\d _\e ^{p+1} \right),
$$
from Lemma \ref{l:53} we deduce that
\begin{align}\label{e:58}
\int _{\O}\d _\e ^{p+1+\e}&= \l _\e ^{\frac{\e
    (n-4)}{2}}c_0^{\e}\left(S^{n/4} - \int _{\R^n
    \diagdown \O }\d _\e ^{p+1}\right)+ O\left(\e \int _{\O}\d _\e
    ^{p+1}Log(1+ \l _\e ^2 |x-a_\e |^2)\right)\nonumber \\
&= (1+o(1))S^{n/4} +O\left((\l _\e d_\e )^{-n}\right) + O(\e ).
\end{align}
Combining \eqref{e:57} and \eqref{e:58}, we obtain the desired result.
\end{pf}\\
\begin{lem}\label{l:55}
Let $u_\e$ satisfying the assumption of Theorem \ref{t:13}. Then we have
$$
{\bf (i)}\quad |u_\e |^{\e}_{L^\infty (\O)} = O(1) \qquad
\quad {\bf (ii)} \quad |v_\e |^{\e}_{L^\infty (\O)} = O(1),
$$
where $v_\e$ is defined in \eqref{e:51}.
\end{lem}
\begin{pf}
We notice that Claim (ii) follows from Claim (i) and Lemma
\ref{l:52}. Then we only need to show that Claim (i) is true.
We define the rescaled functions
\begin{eqnarray}\label{e:59}
\o _\e (y)=M_\e^{-1}u_\e \left(x_\e + M_\e ^{(1-p-\e)/4}y
   \right), \quad y\in \O _\e = M_\e ^{\frac{p-1+\e }{4}}(\O -
x_\e ),
\end{eqnarray}
where $x_\e \in \O$ is such that
\begin{eqnarray}\label{e:510}
M_\e := u_\e (x_\e )= |u_\e |_{L^\infty
  (\O )}.
\end{eqnarray}
$\o _\e$ satisfies
\begin{eqnarray}\label{e:511}
\left\{
\begin{array}{cccccc}
\D^2 \o _\e &=& \o _\e ^{p+\e},& 0< \o _\e \leq 1& \mbox{in}&\O _\e\\
\o _\e (0) &=&1,              &  \D \o = \o _\e = 0    & \mbox{on}&\partial
\O _\e .
\end{array}
\right.
\end{eqnarray}
Following the same argument as in Lemma 2.3 \cite{BEH}, we have
$$
M_\e ^{(p-1+\e )/4}d(x_\e , \partial \O ) \to +\infty \quad
\mbox{as
  } \e \to 0.
$$
Then it follows from  standard elliptic theory
that there exists a positive function $\o $ such that (after passing
to a subsequence) $\o _\e \to \o $ in $C^4_{loc}(\R^n)$, and $\o$
satisfies
\begin{eqnarray*}
\left\{
\begin{array}{ccccc}
\D^2 \o &=&\o ^p,& 0\leq \o \leq 1 & \mbox{in }\R^n\\
\o (0)&=&1,&\n\o (0) =0. &
\end{array}
\right.
\end{eqnarray*}
It follows from \cite{Lin} that $\o$ writes as
$$
\o (y)= \d _{0,\a _n }(y), \quad \mbox{with }\, \a _n =c_0^{2/(4-n)}
$$
Therefore
\begin{eqnarray}\label{e:512}
M_\e ^{\frac{\e (n-4)}{4}}\int _{B(x_\e , M_\e
  ^{\frac{1-p-\e}{4}})}u_\e ^{p+1+\e}(x)dx = \int _{B(0,1)}\o _\e
  ^{p+1+\e}(y)dy \to c>0 \quad \mbox{as } \e \to 0.
\end{eqnarray}
We notice that, as in the proof of Lemma \ref{l:52}
\begin{align}\label{e:513}
\int _{B(x_\e , M_\e ^{\frac{1-p-\e}{4}})}&u_\e ^{p+1+\e}(x)dx
=\int _{B(x_\e , M_\e ^{\frac{1-p-\e}{4}})}(\a _\e P\d _\e
+ v_\e ) ^{p+1+\e}(x)dx\nonumber\\
 & =\a _\e ^{p+1+\e }\int _{B(x_\e , M_\e
  ^{\frac{1-p-\e}{4}})}\d _{a_\e , \l _\e } ^{p+1+\e}(x)dx +
  o(1)\nonumber\\
&=\a _\e ^{p+1+\e } \l _\e ^{\frac{\e (n-4)}{2}}\int _{B(x_\e ,\frac{\l _\e}{ M_\e
  ^{\frac{p-1+\e}{4}}})}\frac{dy}{(1+|y-\l _\e (a_\e -x_\e
  )|^2)^{n+\frac{\e (n-4)}{2}}} + o(1).
\end{align}
Combining \eqref{e:512}, \eqref{e:513} and Lemma \ref{l:52}, we obtain
\begin{eqnarray}\label{e:514}
M_\e ^{\frac{\e (n-4)}{4}}\int _{B(0 ,\l _\e M_\e
  ^{\frac{1-p-\e}{4}})}\frac{dy}{(1+|y-\l _\e (a_\e -x_\e
  )|^2)^{n+\frac{\e (n-4)}{2}}} \to c>0 \quad \mbox{as } \e \to 0.
\end{eqnarray}
We have also
\begin{eqnarray*}
M_\e ^{\frac{\e n}{4}}\int _{B(x_\e , M_\e
  ^{\frac{1-p-\e}{4}})}u_\e ^{p+1}(x)dx = \int _{B(0,1)}\o _\e
  ^{p+1}(y)dy \to c>0 \quad \mbox{as } \e \to 0.
\end{eqnarray*}
Consequently, we find in the same way as above
\begin{eqnarray}\label{e:515}
M_\e ^{\frac{\e n}{4}}\int _{B(0,\l _\e M_\e
  ^{\frac{1-p-\e}{4}})}\frac{dy}{(1+|y-\l _\e (a_\e -x_\e
  )|^2)^{n}} \to c>0 \quad \mbox{as } \e \to 0.
\end{eqnarray}
One of the two following  cases  occurs :\\
{\bf Case 1.}\quad
$\l _\e M_\e ^{\frac{1-p-\e}{4}} \not\to 0 $ as $\e \to 0$.
In this case, we can assume
\begin{eqnarray}\label{e:516}
\l _\e M_\e ^{\frac{1-p-\e}{4}} \geq c_1 > 0, \quad \mbox{as }\e
\to 0.
\end{eqnarray}
The claim follows from \eqref{e:516} and Lemma \ref{l:52}.\\
{\bf Case 2.}\quad
$\l _\e M_\e ^{\frac{1-p-\e}{4}} \to 0 $ as $\e \to 0$. Now
we distinguish two subcases:\\
{\bf Case 2.1.}\quad
$\l _\e |a_\e - x_\e | \not\to +\infty$, as $\e \to 0$.
We can assume that $\l _\e |a_\e -x_\e |$ remains
bounded when $\e \to 0$. Thus, using \eqref{e:515} we obtain
$$
M_\e ^{\frac{\e n}{4}}\left(\l _\e M_\e ^{\frac{1-p-\e
}{4}}\right)^n \to c' > 0 \quad \mbox{as } \e \to 0,
$$
which implies
$$
\l _\e M_\e ^{2/(n-4)} \to c'' > 0 \quad \mbox{as } \e \to 0.
$$
Using Lemma \ref{l:52}, we derive a contradiction : this
subcase cannot happen.\\
{\bf Case 2.2.}\quad
$\l _\e |a_\e - x_\e | \to +\infty$ as $\e \to 0$.
Using \eqref{e:514} and \eqref{e:515}, we obtain
\begin{eqnarray}\label{e:517}
\frac{\l _\e ^n}{\left(\l _\e |a_\e -x_\e |\right)^{2n+\e (n-4)}M_\e
  ^{\e + \frac{2n}{n-4}}} \to C > 0 \quad \mbox{as } \e \to 0
\end{eqnarray}
and
 \begin{eqnarray}\label{e:518}
\frac{\l _\e ^n}{\left(\l _\e |a_\e -x_\e |\right)^{2n}M_\e
  ^{ \frac{2n}{n-4}}} \to C > 0 \quad \mbox{as } \e \to 0.
\end{eqnarray}
 From \eqref{e:517} and \eqref{e:518}, we deduce that
\begin{eqnarray}\label{e:218}
M_\e ^\e \left(\l _\e |a_\e - x_\e |\right)^{\e (n-4)} \to 1
\quad \mbox{as } \e \to 0
\end{eqnarray}
and  we also derive a contradiction. Consequently
Case 2 cannot occur, and the lemma is proved.
\end{pf}\\
Now, arguing as in the proof of Proposition \ref{p:26}, we can easily
derive the following estimate
\begin{lem}\label{l:56}
Let $u_\e$ satisfying the assumption of Theorem \ref{t:13}. Then $v_\e$ occuring
in \eqref{e:51} satisfies
$$
||v_\e ||\leq C \left(\e + \frac{1}{(\l _\e d_\e )^{n-4}}
(\mbox{if }n < 12)+ \frac{Log(\l _\e d_\e )}{(\l _\e d_\e
)^4}(\mbox{if } n=12)+ \frac{1}{(\l _\e d_\e
)^{\frac{n+4}{2}}}(\mbox{if } n > 12)\right)
$$
with $C$ independent of $\e$.
\end{lem}
Now, using the same method in the proof of Claim $(a)$ of Proposition \ref{p:31},
we can easily obtain the following result :
\begin{pro}\label{p:57}
Let $u_\e$ satisfying the assumption of Theorem \ref{t:13}.
Then there exist $C_1>0$ and $C_2>0$
such that
$$
C_1 \frac{H(a_\e ,a_\e )}{\l _\e ^{n-4}}(1+o(1)) + C_2 \e (1+o(1)) =
O\left(\frac{Log(\l _\e d_\e )}{(\l _\e d_\e )^n} +\frac{1}{(\l _\e
    d_\e )^4}  ( \mbox{ if }
  n=5)\right)
$$
where $a_\e$, $\l _\e $ and $d_\e = d(a_\e ,\partial\O )$ are given in \eqref{e:51}.
\end{pro}
We are now able to prove Theorem \ref{t:13}.\\
\begin{pfn}{\bf Theorem \ref{t:13}}
Arguing by contradiction, let us suppose that $(P_{+\e} )$ has a
solution $u_\e$ as stated in Theorem \ref{t:13}. From Proposition
\ref{p:57}, we have
\begin{eqnarray}\label{e:520}
C_1 \frac{H(a_\e ,a_\e )}{\l _\e ^{n-4}} (1+o(1)) +C_2 \e (1+o(1)) =
O\left(\frac{Log(\l _\e d_\e )}{(\l _\e d_\e )^n} + \frac{1}{(\l _\e
    d_\e )^4} (\mbox{ if } n=5
  )\right)
\end{eqnarray}
with $C_1 > 0$ and $C_2 > 0$.\\
Two cases may occur :\\
{\bf Case 1.} $ d_\e \to 0 $ as $\e \to 0$.
Using \eqref{e:520} and the fact that $H(a_\e ,a_\e )
\sim cd_\e ^{4-n}$, we derive a contradiction.\\
{\bf Case 2.}$d_\e \not\to 0 $ as $ \e \to 0$.
We have $H(a_\e ,a_\e ) \geq c >0 $ as $\e \to 0$ and
\eqref{e:520} also leads to a contradiction. Thus our result follows.
\end{pfn}

\end{document}